\documentclass{amsart}
\usepackage{epsfig}
\usepackage{amssymb}
\begin{document}

\def\square{\hfill ${\vcenter{\vbox{\hrule height.4pt
	    \hbox{\vrule width.4pt height7pt \kern7pt
	        \vrule width.4pt}
	    \hrule height.4pt}}}$}

\def\Re{{\text{Re}}}
\def\Im{{\text{Im}}}

\newenvironment{pf}{\noindent {\it Proof:}\quad}{\square \vskip 12pt}

\newtheorem{thm}{Theorem}
\newtheorem{cor}[thm]{Corollary}

\title{A Note on Strong Geometric Isolation in $3$--Orbifolds}
\author{Danny Calegari}
\maketitle

\begin{abstract}
Neumann and Reid in \cite{NeRe}
describe a $2$--cusped hyperbolic $3$--orbifold in which the cusps
are geometrically isolated. Based on numerical evidence provided by Jeff Weeks'
{\tt snappea} 
program, they conjecture that the cusps are strongly geometrically
isolated, a fact which we establish here. We also give a 
parameterization of the
Dehn Surgery Space of this orbifold which has amusing properties.
\end{abstract}

\section{Geometric Isolation}

Following \cite{NeRe} we make the following definitions:

Given a hyperbolic $3$-orbifold $M$ with cusps $c_1, \dots c_k, c_{k+1} \dots
c_h$ we say that cusps $1, \dots, k$ are {\it geometrically isolated} from
cusps $k+1, \dots ,h$ if any deformation of the hyperbolic structure on $M$ 
induced by Dehn filling cusps $k+1, \dots, h$ while keeping cusps $1,\dots,k$
complete does not change the Euclidean structure at cusps $1, \dots, k$.

We say that cusps $1, \dots,k$ are {\it first order isolated} from cusps
$k+1,\dots,h$ if the map from the space of deformations induced by Dehn
filling cusps $k+1, \dots, h$ while keeping cusps $1,\dots,k$ complete to
the space of Euclidean structures at cusps $1,\dots, k$ has zero derivative
at the point corresponding to the complete structure.

We say that cusps $1, \dots,k$ are {\it strongly geometrically isolated}
from cusps $k+1,\dots,h$ if for any fixed Dehn filling on cusps
$1, \dots,k$, the geometry of the (possibly filled) cusps $1, \dots,k$ is
unchanged by any Dehn filling on cusps $k+1,\dots,h$.

It is immediate from these definitions that strong isolation implies
geometric isolation, which in turn implies first order isolation.

In \cite{NeRe} it is shown that first order isolation and strong isolation
are symmetric conditions in the sets $1,\dots,k$ and $k+1,\dots,h$, and
that they can be given an analytic definition in terms of the $\Phi$ function
defined in \cite{NeZa}.

\section{The orbifold $A$ and the manifold $A^*$}

Neumann and Reid in \cite{NeRe}
construct a $2$-cusped hyperbolic orbifold $A$ with geometrically
isolated cusps. Based on numerical evidence, they conjectured that
the cusps are strongly isolated. The orbifold $A$ is particularly interesting
for two reasons: firstly it is arithmetic, and secondly, all other known
examples of manifolds with strongly isolated cusps contain a rigid totally
geodesic separating surface which ``forces" the isolation. There is
evidence that this example does not contain such a surface, based on the fact
that its volume is so small. In this paper we find explicit formulae describing how
the tetrahedra making up an ideal triangulation of $A$ are distorted as
the hyperbolic structure on $A$ is deformed which suggests that the existence of
such a surface is unlikely.

$A^*$ is a $4$--cusp link complement in $S^3$. $A$ is obtained from $A^*$ by
$(2,0)$ Dehn surgery on two of the links. It is a curious fact, which we will
exploit in our calculations, that $A^*$ is a double cover of $A$, obtained by
a twofold branching of $S^3$ over the two filled cusps.

\begin{figure}
\scalebox{.35}{\includegraphics{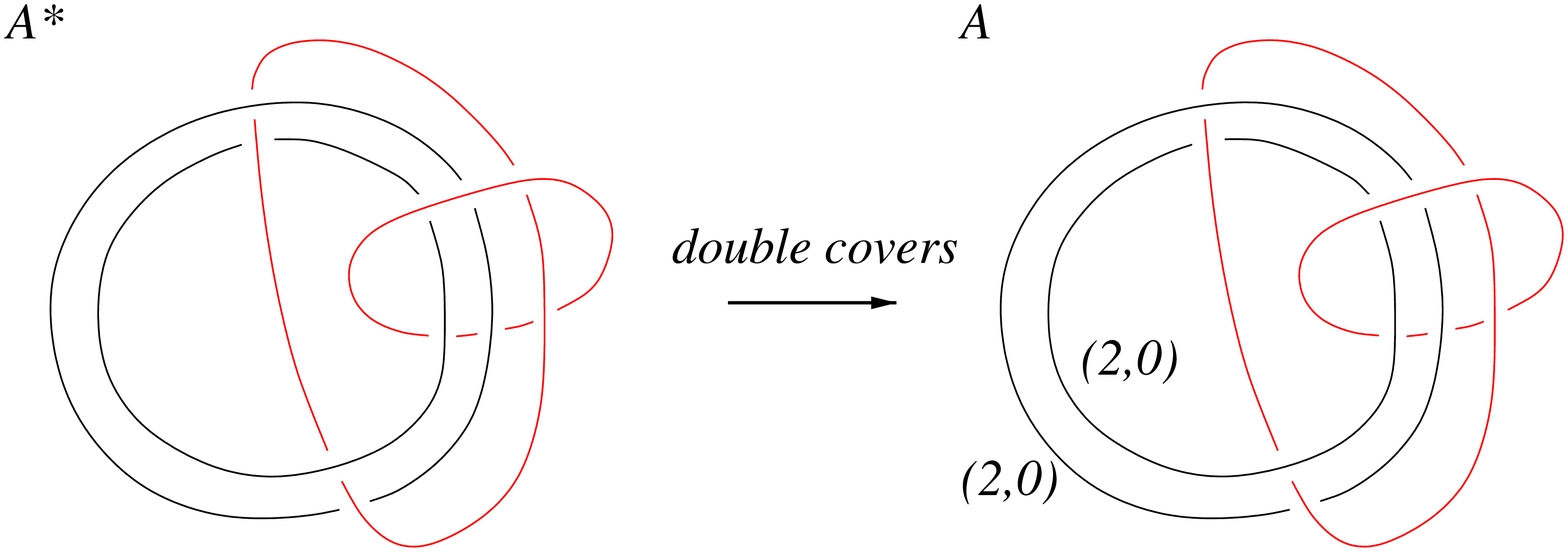}}
\caption{The orbifold $A$ and the manifold $A^*$}
\end{figure}

We can draw $A^*$ as a $T^2 \times (0,1)$ complement. The two removed cusps are
linked as in Figure 2. The two cusps in the figure we label as $X$ and $Y$. $W$ 
is the cusp $T^2 \times \lbrace 0 \rbrace$, and $Z$ the cusp 
$T^2 \times \lbrace 1
\rbrace$. As in \cite{NeZa} we let $u_i$ and $v_i$ denote 
the logarithm of the
holonomy of $m_i$ and $l_i$, the meridian and longitude 
respectively of cusp $i$.
Rotating through an angle $\pi/2$ gives the same
link complement after isotopy. 

In fact, we can see that $A^*$ retains its $4$--fold symmetry even after
cusps $X$ and $Y$ are $(m,n)$--filled. Hence the cusp shapes at 
$T^2 \times \lbrace 0 \rbrace$ and $T^2 \times \lbrace 1 \rbrace$ 
are constant. By taking the
quotient we see that in $A$ the two cusps are geometrically isolated. Note that
this technique can be generalized. Let $M$ be any $4$--fold symmetric
$T^2 \times (0,1)$ 2--link
complement for which there is an orientation--preserving isometry permuting each
pair of cusps. Quotienting by this isometry we obtain an orbifold
in which one cusp is geometrically isolated from the other. Alternatively, let
$M$ be any $4$--fold symmetric $T^2 \times (0,1)$ knot complement. Then the
cusp corresponding to this knot is geometrically isolated from the other two
cusps. It is not difficult to construct infinitely many examples this way.

\begin{figure}
\scalebox{.45}{\includegraphics{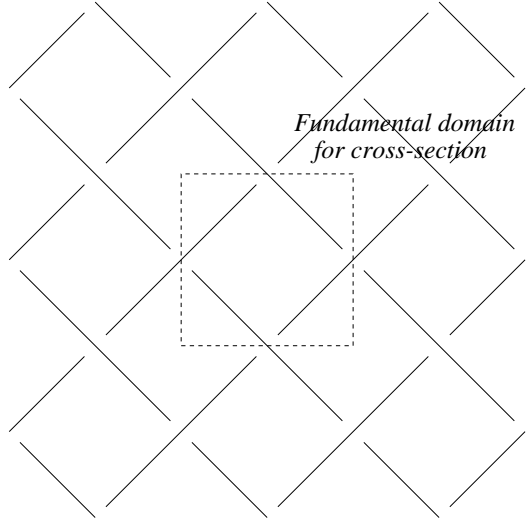}}
\caption{$A^*$ as a $T^2 \times [0,1]$ complement}
\end{figure}

\section{Strong Isolation of Cusps}

$A^*$ has a hyperbolic structure with fundamental domain given by two regular
ideal hyperbolic octahedra identified along their boundaries. 
To see this, note that
$A^*$ consists of an alternating link complement in $T^2 \times (0,1)$, 
and therefore
by general properties of alternating links, 
a fundamental domain for $A^*$ can be
found by two--coloring the surface on which the link 
projection lies (in this case a
torus), and gluing the regions above and below this surface with a twist of
$\pi/2$ or $-\pi/2$ according to the coloring (see \cite{AiRu}).
The link in $A^*$ separates the
dividing torus into two squares. 
The suspension on each of these squares for which
the suspension points are cusps $W$ and $Z$ is an octahedron,
from which the vertices are removed, since they lie on the cusps. 
Gluing the top four
faces of each octahedron in the obvious manner, 
the bottom four faces of each are
identified with a twist of $\pi$, relative to the 
top identification. The pattern of
identifications is given in Figure 3.

To prove that $A$ is strongly geometrically isolated, it suffices to show that
$v_1(u_1,u_3)$ and $v_3(u_1,u_3)$ are functions
solely of $u_1$ and of $u_3$ respectively; i.e. that we can define
$v_1(u_1)$ and $v_3(u_3)$.

We subdivide each octahedron into $4$ ideal tetrahedra by cutting along the
two faces $XZXW$ and $YZYW$, 
so that each tetrahedron has vertices $X, Y, Z$ and
$W$. We denote the simplex parameters of the four tetrahedra from the
first octahedron as $z_1,z_2,z_3,z_4$ relative to the central edge, and
for the tetrahedra from the second octahedron as $w_1,w_2,w_3,w_4$.
The triangulated cusp shapes for the complete structure are given in Figure 4.

In the notation of \cite{NeZa} we have the consistency relations

$$w_1w_2w_3w_4=1, \; z_1z_2z_3z_4=1$$
$$w_1''w_2'w_3''w_4'z_1'z_2''z_3'z_4''=1, \; w_1w_3z_1z_3=1$$

The holonomies of $l$ and $m$ for the four cusps are

$$l_W:w_4''{z_2''}^{-1}{z_1'}^{-1}w_1', \; m_W:z_2'{w_4'}^{-1}{w_3''}^{-1}z_3''$$
$$l_Z:w_2''{z_2''}^{-1}{z_3'}^{-1}w_1', \; m_Z:z_2'{w_2'}^{-1}{w_3''}^{-1}z_1''$$
$$l_Y:z_3''{z_2''}^{-1}{w_2'}^{-1}w_3', \; m_Y:z_2'{z_3'}^{-1}{w_1''}^{-1}w_4''$$
$$l_X:z_3''{z_4''}^{-1}{w_2'}^{-1}w_1', \; m_X:z_4'{z_3'}^{-1}{w_3''}^{-1}w_4''$$

\begin{figure}
\scalebox{.35}{\includegraphics{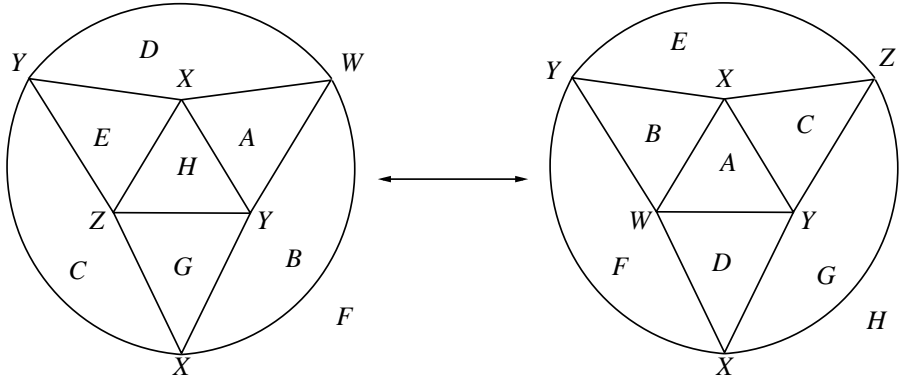}}
\caption{Identifying the faces of two regular ideal octahedra in this
pattern gives $A^*$}
\end{figure}

Since we are lifting Dehn surgeries from $A$, we require that

$$u_W=u_Z, \; u_Y=u_X$$

from which, and from the Dehn surgery equations it will follow that

$$v_W=v_Z, \; v_Y=v_X$$

As in \cite{NeZa} we denote $z' = \frac{z-1} z$ and $z'' = \frac 1 {1-z}$.

Reparameterizing, we find that we can write everything in terms of two complex
variables, $\alpha$ and $\beta$. We can define the simplex parameters in terms
of these two variables as follows:

$$w_4 = \frac {\alpha i + \alpha + 1 - i} {\alpha i - \alpha + 1 - i}, \;
z_1 =  \frac {\alpha + \alpha i} {2 - \alpha + \alpha i}$$
$$w_2 = \frac {\alpha i + \alpha -1 - i} {\alpha i - \alpha +1 + i}, \;
z_3 = \frac {\alpha i + \alpha - 2i} {\alpha i - \alpha}$$
$$w_1 = \frac {2i - \beta - \beta i} {\beta - \beta i}, \;
z_2 = \frac {i - 1 - \beta - \beta i} {\beta - \beta i + i - 1}$$
$$w_3 = \frac {- \beta - \beta i} {\beta - \beta i - 2}, \;
z_4 = \frac {1 + i - \beta - \beta i} {\beta - \beta i - i - 1}$$

With these definitions, 
it can easily be verified that $l_W=l_Z$ and $m_W=m_Z$ are
functions only of $\beta$, and that $l_Y=l_X$ and $m_Y=m_X$ 
are functions only of
$\alpha$; for example, 
$l_W= {z_2''}^{-1}w_1'w_4''{z_1'}^{-1}=\frac {w_1'} {2z_2''}$, which
is a function only of $\beta$. Explicitly we can calculate

$$l_W = \frac {\beta(\beta-i)} {(\beta-1)(\beta - (1+i))}, \;
m_W = \frac {\beta(\beta-1)} {(\beta-i)(\beta - (1+i))}$$
$$l_Y = \frac {\alpha(\alpha-1)} {(\alpha - i)(\alpha-(1+i))}, \;
m_Y = \frac {(\alpha-1)(\alpha - (1+i))} {\alpha(\alpha-i)}$$

We therefore have the following:

\begin{thm}
The orbifold $A$ has strongly isolated cusps.
\end{thm}

\begin{figure}
\scalebox{.5}{\includegraphics{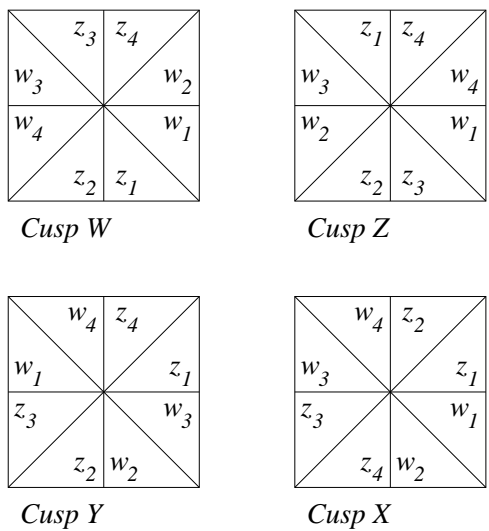}}
\caption{Triangulated cusp shapes for the complete structure}
\end{figure}

\section{Geometric Description of Reparameterization}

There is a nice geometric definition of the 
parameterization in terms of the two
complex numbers $\alpha$ and $\beta$, which we now describe.

Given a square $ABCD$ in the complex plane of side length $1$, with vertices at
$0,1,1+i,i$, we choose a point $O$ at the complex number $\alpha$. We construct
the points $R$, $S$, $T$ and $U$ so that the triangles 
$ASO$, $BTO$, $CUO$ and $DRO$
are similar with angles $\frac \pi 4, \frac \pi 2$ and $\frac \pi 4$, 
as illustrated in Figure 5.
Then it is an easy calculation to see that as complex numbers,
$T = R+1$ and $S=U+i$. 
Therefore the octagon $ASBTCUDR$ tiles the complex plane.

The four triangles $DUO$, $CTO$, $BSO$ and $ARO$ are 
precisely (up to similarity) the
horoball sections of the simplices 
$z_2$, $w_3$, $z_4$ and $w_1$ respectively after
Dehn filling cusps $Z$ and $W$.

There is a similar, though ``mirror reversed'' picture 
describing how $w_2$, $w_4$, $z_1$ 
and $z_3$ are determined as a function of $\alpha$. From the fact that this
octagon tiles the plane, the consistency relations $w_1w_3=z_2z_4=-1$ and
$w_1''w_3''z_2''z_4''=-\frac 1 4$ are geometrically obvious. 
From the mirror-reversed
picture, the consistency relations $w_2w_4=z_1z_3=-1$ and $w_2'w_4'z_1'z_3'=-4$
can similarly be deduced.

A relation of the form $w_4''{z_1'}^{-1}=\frac 1 2$ 
is equivalent to the geometric fact
that the orientation preserving
similarity taking $AO$ to $OS$ takes $R$ to the midpoint of $OB$. From simple
relations such as these, 
the independence of $l_W$ and $m_W$ from $\alpha$ and of
$l_Y$ and $m_Y$ from $\beta$ can be seen.

\begin{figure}
\scalebox{.45}{\includegraphics{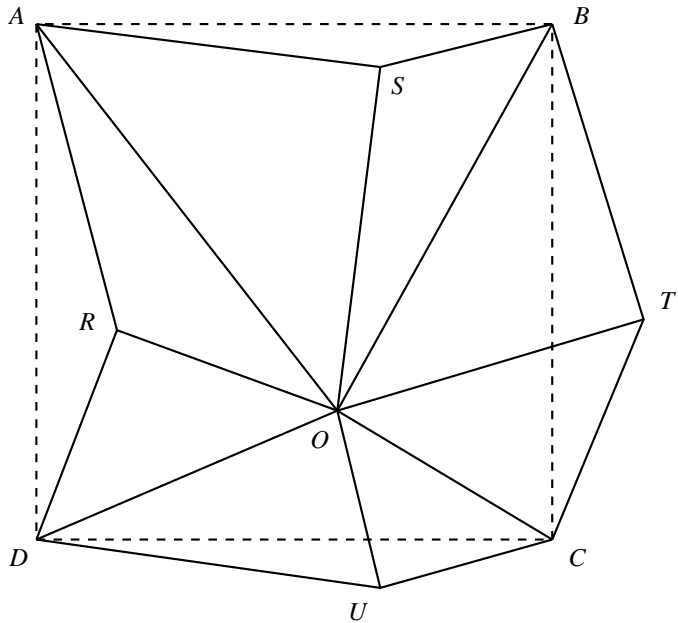}}
\caption{Geometric description of reparameterization}
\end{figure}

\section{The Dehn Surgery Parameter Space}

There is an interesting relationship between the 
parameterization of Dehn Surgery Parameter 
Space by real pairs $p_i,q_i: i=1,2$, and by $\alpha, \beta$. 
In particular, we have the
following

\begin{thm}\label{circle_par}
The circle $$|\beta- 1/2 - i/2| = \frac 1 {\sqrt 2}$$
corresponds to the values of
$(p,q)$ lying on the square with vertices at $(\pm 2,\pm 2)$.
\end{thm}
\begin{pf}
The circle $|\beta - 1/2 - i/2| = 1/\sqrt 2$ 
passes through the square in
the complex plane with vertices at $0,1,i,1+i$. 
Therefore, by the description in
section \S 4, it corresponds precisely to the values of $\beta$ 
for which four of the
simplices are degenerate. 
That is to say, the simplex parameters $w_1,w_3,z_2,z_4$ are
real. By symmetry, we need only consider the 
case that $\beta$ lies on a fixed ``quadrant''
of the circle, for the sake of argument, the arc 
between $D$ and $C$ in Figure 5.

However we know that
$$v_W=\log(1-z_2) + \log(w_1-1) - \log(w_1) - \log(2)$$
$$u_W=\log(z_2-1) + \log(1-w_3) - \log(z_2) - \log(2)$$
and moreover, that
$$p_2 u_W + q_2 v_W = 2 \pi i$$

But then, by inspecting Figure 5 it can be seen that 
for $\beta$ taking values on this
arc, $w_1 \le -1$, $z_2>1$ and $0<w_3<1$, so that  $\Im(v_W), \Im(u_W)$
must take the values of $\pi, 0$ respectively. Hence, $q_2 = 2$. 
It is a simple matter to check that $|\Re(v_W)| \le |\Re(u_W)|$, 
and therefore that $|p_2| \le 2$.

By a continuity argument,
it can be seen that vertices $C$ and $D$ correspond to the points
$(-2,2)$ and $(2,2)$ in $p_2,q_2$ space, and the proof follows.
\end{pf}
Notice of course that a corresponding fact holds for $p_1,q_1$ and $\alpha$.

\begin{cor}
$Vol(A_{(p_2,q_2)})=Vol(A)/2$ for $(p_2,q_2)$ taking values on the square
with vertices as in Theorem~\ref{circle_par}.
\end{cor}

By a similar argument we can show

\begin{thm}
The ``circle'' $|\beta| = \infty$ corresponds to the values of $(p_2,q_2)$
lying on the square with vertices $(\pm 1,\pm 1)$
\end{thm}
\begin{pf}
For the sake of argument, the infinite circle can 
be understood as a limit of larger and
larger circles, in the Hausdorff topology on the 
Riemann sphere $\widehat{\Bbb C}$.

It is easy to see that as $\beta \to \infty$, 
the simplex parameters
$$z_2, z_4, w_1, w_3 \to -i$$ 
(that is to say, their orientation and volume is negative).
This implies that each of $v_W, u_W$ approaches 
the value of $0$ or $\pm 2 \pi i$ according
to which branch of the logarithm is taken. 
This depends on which homotopy path $\beta$
takes from $1/2 + i/2$ to $\infty$ in the 
$4$-punctured sphere $\widehat{\Bbb C} - \lbrace
0,1,i,1+i \rbrace$. We restrict the values of $\beta$ for the moment
to the simply connected domain
consisting of ${\Bbb C}$ minus the four infinite rays 
emanating from the four corners
of the square. The appropriate branch of the 
logarithm can then be determined by
observing how the simplex parameters deform as 
$\beta \to \infty$ and then requiring that
$u_W, v_W$ be continuous. It remains to check that $|p_2|,|q_2| \le 1$ in the
limit. Without loss of generality, we choose $\beta$ with
$ \Re(\beta) > \Im(\beta) \gg 0 $. Then $p_2 \to 1$, 
by the argument above. It suffices
to check that $|\Re(v_W)| \ge |\Re(u_W)|$ as $\beta \to \infty$. But by the
parameterization in terms of $\beta$, 
it is immediate that $|l_W|>|m_W|>1$. The case
for other large absolute values of $\beta$ follows by symmetry, 
and a continuity argument for the cases that 
$|\Re(\beta)| = |\Im(\beta)|$. The proof follows.
\end{pf}

One must be wary about assuming that geometric structures corresponding to
the case that $(p_2,q_2)$ lie in the interior of the square with vertices as
in Theorem~\ref{circle_par} exist.

For the case that $(p_2,q_2)$ lie on this square, four of
the tetrahedra are flat and disjoint except at the vertices. The other four
tetrahedra can be glued to the flat four, to produce an incomplete hyperbolic
manifold, by the consistency equations. 
The fact that each flat tetrahedra has a
(geometric) bicollared neighborhood in the manifold 
so obtained not intersecting
any of the other four tetrahedra implies that one can ``pump some air'' 
into each
of the four flat tetrahedra without changing the homeomorphism type of the
manifold. That is to say, the incomplete manifold 
so obtained is homeomorphic to
$A^*$, and the complete manifold obtained by 
performing geometric Dehn surgery is
homeomorphic to a manifold obtained by 
topological Dehn surgery on $A^*$. Note by
a theorem proved in \cite{H} one can find 
geometric structures corresponding to
$(p_2,q_2)$ in an open neighborhood of
$\lbrace \text{square } - \text{ vertices}\rbrace$ --- that is to say,
including some cases in which four of the tetrahedra are negatively oriented.

Notice that as $\beta$ approaches one of $0,1,i,1+i$ that 
one of the $(1,1,\sqrt2)$ triangles
in a horoball section gets arbitrarily small 
with respect to the other three. This is to
say, the length of the filled geodesic in a 
complete structure (if one were to exist) goes
to infinity. This suggests that the orbifolds 
$A_{(\pm 2,\pm 2)}$ contain incompressible
tori or annuli, a fact which can be verified geometrically.

\section{Acknowledgements}

The author wishes to express his gratitude to Walter Neumann and Iain Aitchison
with whom he had many very fruitful conversations regarding 
the material in this paper. 
He particularly wants to thank Walter Neumann for his generosity in taking
the time to make numerous useful remarks on the presentation 
and details of this material.

\end{document}